\documentclass[draft,leqno,12pt]{article}

\usepackage{amsmath}
\usepackage{amscd}
\usepackage{amssymb}
\usepackage{enumerate}
\usepackage{indentfirst}
\usepackage{latexsym}
\usepackage{multicol}
\usepackage{theorem}
\usepackage{esint}

\renewcommand{\rho}{\varrho}

\renewcommand{\div}{\operatorname{div}}
\newcommand{\grad}{\nabla}

\makeatletter\@addtoreset{equation}{section}\makeatother

{\theorembodyfont{\rmfamily}
\newtheorem{remark}{Remark}[section]}\newtheorem{theorem}{Theorem}[section]
\newtheorem{proposition}{Proposition}[section]
\newenvironment{proof}{\textit{Proof. }}{\hfill$\Box$}
\newenvironment{proofoftheorem}[1]{\noindent{\bf Proof of Theorem #1.
}}{\hfill$\Box$}

\newcommand{\LlnL}{\text{\bf{L}}\ln\text{\bf{L}} }

\begin{document}


\renewcommand{\thefootnote}{\fnsymbol{footnote}}

\title{Zygmund spaces, inviscid limit and uniqueness of Euler flows
\footnotetext{\textbf{Mathematics Subject Classification (2000). }
{76B03, 76D09, 35Q30.}\hfill\break}
\footnotetext{\textbf{Keywords. }{Zygmund spaces, Euler equations, uniqueness, Navier-Stokes equations, inviscid limits.}\hfill\break}}
\author{\textsc{P.B. Mucha$^1$} and \textsc{W.M. Rusin$^2$}
}
\date{}
\maketitle

\begin{center}
{\small {

1. Institute of Applied Mathematics and Mechanics,Warsaw University}

{ul. Banacha 2, 02-097 Warszawa, Poland }

{E-mail: {\tt p.mucha@mimuw.edu.pl}}

\medskip

{2. School of Mathematics, University of Minnesota}

{206 Church Street SE, Minneapolis, 55455 MN, USA }

{E-mail: {\tt rusin018@math.umn.edu}}

}

\end{center}

\begin{abstract}

The paper improves the classical uniqueness result for the Euler system in the $n$ dimensional case
assuming that 
$\nabla u^E \in L_1(0,T;BMO(\Omega))$, only. Moreover the rate of the convergence for the 
inviscid limit of solutions to the Navier-Stokes equations is obtained, provided the same
regularity of the limit Eulerian flow. A key element of the proof is a logarithmic inequality 
between the Hardy and $L_1$ spaces which is a consequence of the basic properties of the Zygmund space
$\LlnL$.

\end{abstract}



\section{Introduction}

The analysis of the evolutionary Euler system modeling the motion of incompressible flows in $n$ dimensional bounded domains
 is the subject of this paper. We want to study the issue of uniqueness and the problem of the inviscid limit for the 
 Navier-Stokes equations treated as an approximation of the system of inviscid flows.
 
 The classical results \cite{KAT1} and \cite{Y1} require that solutions to the Euler system should belong at least to the class of
 regularity which guarantees that the velocity is in the class $u^E \in L_1(0,T;W^1_\infty(\Omega))$. Due to that fact we  obtain the 
 following estimate
 \begin{equation}\label{1.1}
 \left|\int_0^T\int_\Omega v\cdot \nabla u^E v dxdt\right| \leq 
 C\|\nabla u^E\|_{L_1(0,T;L_\infty(\Omega))}\|v\|^2_{L_\infty(0,T;L_2(\Omega))}
 \end{equation}
 which is the core of  methods in \cite{KAT1},\cite{Y1}. Having inequality (\ref{1.1}) the uniqueness of solutions to 
 the Euler system follows from elementary energy estimates.
 
 The goal of our paper is to improve the classical approach to the Euler system 
 replacing $L_\infty$ by the $BMO$ space. Because of relatively low 
 regularity in the studied problem we cannot apply the properties of the $BMO$ space directly. A key element of our technique will be 
 an application of  properties of Zygmund spaces $\LlnL$ (see \cite{Zyg1}). This analysis enables us to prove the following bound 
 \begin{equation}\label{1.2}
 \|w\|_{ {\mathcal H}^1(\Omega)}\leq C\|w\|_{L_1(\Omega)}\left[
 |\ln \|w\|_{L_1(\Omega)}| + \ln (1+\|w\|_{L_\infty(\Omega)})\right]
 \end{equation}
which
 measures the difference between the $L_1$ and Hardy space ${\mathcal H}^1$.

 We will study
 the uniqueness criteria which play an important role in analysis based on weak solutions where -- by the definition --  the high 
 regularity is not admitted. They allow consideration of larger class of  external (initial) data to obtain information almost 
 the same as for smooth data. Thanks to (\ref{1.2}) we will be able to prove that the criteria for the incompressible 
 Euler system should guarantee that $\nabla u^E \in L_1(0,T;BMO(\Omega))$, only, replacing the stronger condition from (\ref{1.1}).
  Moreover the analysis will enable to consider
 the approximation of solutions to the Euler system by solutions to the Navier-Stokes equation with small viscosity coefficient.
 The $BMO$ space is specially distinguished in two dimensions, since it is the limit space for the imbedding
 $H^1({\mathbb R}^2)\subset BMO(\mathbb R^2)$   (we can not obtain $L_\infty$ here). This case is of our special interest, since
 we are able to point out good examples for which the $L_\infty$ regularity with respect to spatial coordinates is too strong.

 Our first result, being the fundamental tool of analysis of the Euler system, is the following
 
 \smallskip
 
 \begin{theorem}\label{estimate_i}
	Let $f \in BMO(\mathbb{R}^n)$ and $g \in L_1(\mathbb{R}^n) \cap L_\infty(\mathbb{R}^n)$, then
\begin{equation}\label{ineq_est}
	\left|\int_{\mathbb{R}^n}fg\;dx \right|\leq C\|f\|_{BMO(\mathbb{R}^n)}\|g\|_{L_1(\mathbb{R}^n)}
	\left[  |\ln\|g\|_{L_1(\mathbb{R}^n)} |+\ln(1+\|g\|_{L_\infty(\mathbb{R}^n)})
	  \right].
\end{equation}
\end{theorem}
 
Theorem 1.1 is a version of the logarithmic Sobolev inequality for the Hardy  and $L_1$ spaces.
The structure of (\ref{ineq_est}) and its proof is essentially based on the properties of the Zygmund spaces
$\LlnL$ -- see \cite{TOR1},\cite{Zyg1}. Inequality (\ref{1.2}) (or (\ref{ineq_est})) can be compared with a similar estimate between 
the $L_\infty$ and $BMO$ spaces from \cite{KOZ}. The authors have shown that
\begin{equation}\label{1.4}
\|w\|_{L_\infty(\Omega)}\leq C\left[1+\|f\|_{BMO(\Omega)}(1+\ln^{+} \|f\|_{W^s_p(\Omega)})\right] \mbox{ \ \ for \ }
s>\frac np.
\end{equation}
Inequality (\ref{1.4}) helped to improve the classical result for the Euler system for 
the blow-up criteria replacing $L_\infty$ by the $BMO$ space. In our problems we can not apply estimate (\ref{1.4}),
but it stays the motivation for Theorem 1.1.

We want to apply Theorem \ref{estimate_i} to analyze the Euler system
\begin{equation}\label{Euler}
	\begin{gathered}
	u^E_t+u^E\cdot\nabla u^E +\nabla p^E = 0 \hfill \text{ \ \ \ in \ }  \Omega\times(0,T), \\
	\text{div }u^E=0 \hfill \text{ in \ }  \Omega\times(0,T), \\
	\vec n\cdot u^E =0 \hfill \text{ on \ }  \partial\Omega \times (0,T), \\
	u^E|_{t=0}=u_0 \hfill \text{ in \ }  \Omega.
	\end{gathered}	
	\end{equation}
where $u^E:\Omega \times (0,T) \to \mathbb R^n$ is the velocity field, $p^E:\Omega \times (0,T) \to \mathbb R$ is the pressure, 
$\vec n$ - unit, outward normal to $\partial\Omega$,
 $u_0:\Omega \to \mathbb R^n$ - a divergence-free initial velocity field. We exclude the external forces 
 (the r.h.s. of $(\ref{Euler})_1$), since it would not
 provide any new analytical difficulties, but only few technical elementary estimates.
 
We prove the following result concerning the issue of uniqueness of solutions to the system (\ref{Euler})

\smallskip
 
\begin{theorem}\label{uniqueness_i} Let $\Omega$ be bounded domain in $\mathbb R^n$ with smooth boundary.
Let $u^E_1$ and $u^E_2$ be two solutions to the Euler system (\ref{Euler}) with initial data $u_0$ such that
\begin{equation}\label{1.6}
\begin{array}{c}
\nabla\;  u^E_1, \nabla\; u^E_2 \in L_1(0,T;BMO(\Omega)) \mbox{ \ \ and \ \ }\\[8pt]
 u^E_1,u^E_2 \in L_\infty(0,T;L_{2+\sigma}(\Omega)) \mbox{ \ for given \ }
\sigma >0,
\end{array}
\end{equation}
 then $u^E_1\equiv u^E_2$.
\end{theorem} 
 
 \smallskip
 
The above result generalizes the classical theory. The proof of Theorem 1.2 is essentially based on Theorem 1.1 and
an application of the Osgood theorem (see \cite{HART}) giving uniqueness in the ODE theory. 
Additionally, the properties of the $BMO$-space allow to exchange 
 the gradient $\nabla u$  by the vorticity $rot\;u$ in the condition (\ref{1.6}) which may simplify the application of Theorem 1.1.

An improvement of the classical results as in \cite{KAT1}  and \cite{Y1} can be found \cite{Vish} and \cite{Y1}. 
However the relaxation
of the $L_\infty$-regularity in the space (in \cite{Vish} we even have a bit weaker space than $BMO$) forces that the regularity with respect to time is required  to belong to the $L_\infty$-class.
 So the regularity with respect to time has to be even
stronger than in (\ref{1.1}). Our approach enables to
keep the weak condition in the $L_1$-norm in $(0,T)$. Moreover thanks to Theorem 1.1 we omit numerous technical estimates 
which often appear in results of that type.
 
 \smallskip 
 
Our last result concerns the inviscid limit for solutions to the Navier-Stokes equations under
slip boundary conditions 
\begin{equation}\label{NS}
	\begin{gathered}
    	u^\nu_t+u^\nu\cdot\grad u^\nu - \div \mathbb{T}(u^\nu,p^\nu) = 0 \hfill \text{ \ \ \ \ in } \Omega\times (0,T),\\
    	\div u^\nu=0 \hfill \text{ in } \Omega\times (0,T), \\
	\vec n\cdot \,\mathbb{T}(u^\nu,p)\cdot \,\vec\tau + \alpha\, u^\nu\cdot \vec \tau =0 \hfill \text{ on }\partial\Omega \times (0,T), \\
    	n\cdot u^\nu =0 \hfill \text{ on } \partial\Omega \times (0,T), \\
    	u^\nu|_{t=0}=u_0 \hfill \text{ in } \Omega,
	\end{gathered}
	\end{equation}
where $u^\nu:\Omega \times (0,T) \to \mathbb R^n$ is the velocity field, $p^\nu:\Omega \times (0,T) \to \mathbb R$ is the pressure,
 $\vec n$ - unit, outward normal to $\partial\Omega$,
 $\vec \tau$ - unit, tangent to $\partial\Omega$, $\mathbb{T}(u^\nu,p^\nu)=\nu\mathbb{D}(u)-p^\nu\,\text{Id}$ - stress tensor, 
 $\alpha$ - describes friction coefficient of the boundary, $u_0:\Omega \to \mathbb R$ - a divergence-free initial velocity field. 
 We are able to consider 
 different boundary conditions than $(\ref{NS})_{3,4}$, however by the results from \cite{Cl},\cite{MU}
  the form of (\ref{NS}) seems to be the most suitable for the issue.
 
We prove

\smallskip 

\begin{theorem}Let $\Omega$ be bounded domain in $\mathbb R^n$ with smooth boundary,
	let $u^\nu$ be a solution to the Navier-Stokes system (\ref{NS}) and $u^E$ be the solution to the Euler system (\ref{Euler})
	both with initial data $u_0$.
	Fix $T>0$, $\sigma>0$ and consider $\nu$ such that $0< \nu \leq \nu_0$. Assume that
\begin{equation}\label{1.8}
\begin{array}{c}
\|\nabla u^E(\cdot,t)\|_{BMO(\Omega)}\leq f_0(t) \mbox{ \ \ and \ \ } f_0\in L_1(0,T),\\[10pt]
\|\nabla u^\nu(\cdot,t)\|_{L_2(\Omega)}+\|\nabla u^E(\cdot,t)\|_{L_2(\Omega)}\leq g_0(t) \mbox{ \ \  and \ \ } g_0 \in L_2(0,T),\\[10pt]
\|u^\nu(\cdot,t)\|_{L_{2+\sigma}(\Omega)}+\| u^E(\cdot,t)\|_{L_{2+\sigma}(\Omega)}\leq h_0(t) 
 \mbox{ \ \  and \ \ } h_0 \in L_\infty(0,T).
\end{array}
\end{equation}	
Then considering the inviscid limit of solutions to (\ref{NS}) we obtain
\begin{equation}\label{1.9}
\sup_{0\leq t\leq T} \|(u^\nu-u^E)(\cdot,t)\|_{L_2(\Omega)} \to 0 \mbox{ \ \ \ as \ \ \ } \nu \to 0^+,
\end{equation}
where the precise rate can be expressed by the properties of functions $f_0,g_0$ and $h_0$.

Additionally if we assume extra that
\begin{equation}\label{1.10}
\sup_{0\leq t \leq T} |f_0(t)+g_0^2(t)| \leq M,
\end{equation}
then we obtain the following explicit rate of convergence
\begin{equation}\label{1.11}
\sup_{0\leq t \leq T}\|(u^\nu-u^E)(\cdot,t)\|_{L_2(\Omega)}  \leq C\nu^{e^{-2MT}}.
\end{equation}

\end{theorem}
 
\smallskip

Theorem 1.3 gives general conditions for the inviscid limit to solutions of the Navier-Stokes equations, provided very low
(lowest known) conditions on the regularity of solutions to (\ref{NS}) 
with respect to the viscous coefficient $\nu$. The main disadvantage is that in the general case we  are not able to construct
solutions fulfilling (\ref{1.8}). However in a special case in two dimensions (see \cite{Ru} and \cite{MU} for the case with homogeneous
boundary data) we find a class of solutions to (\ref{NS}) that fulfills assumptions (\ref{1.10}).
 Then by Theorem 1.3 we obtain an explicit rate of convergence to solution of
the Euler system given by (\ref{1.11}). A similar result has been known only for the whole two dimensional case \cite{CONST} under the 
classical assumption
$\nabla u^E \in L_1(0,T;L_\infty(\mathbb R^2))$. The the rate of convergence of $\sup_{0\leq t \leq T}\|u^\nu-u^E \|_{L_2(\mathbb{R}^2)}$ is estimated by $\sim \sqrt{\nu T}$, however the initial data considered in \cite{CONST} correspond to a vortex patch -- vorticity is localized to a bounded domain with smooth boundary. 

Maybe there is a hope to find  a realization  of (\ref{1.8}) by some class of  solutions to the Navier-Stokes equations.
However the problem seems to be challenging.

In the proceeding, $\Omega$ will be always understood as a bounded subset of $\mathbb{R}^n$ with smooth boundary $\partial \Omega$. 
Spaces $(L_p(\Omega),\|\cdot\|_{L_p(\Omega)}),(L_p(\mathbb{R}^n),\|\cdot\|_{L_p(\mathbb{R}^n)})$ for $p\in [1,\infty]$ denote the 
usual Lebesgue spaces. Spaces $BMO(\mathbb{R}^n)$ and $BMO(\Omega)$ are understood as spaces of measurable functions for which
 corresponding semi-norms
$$
||f||_{BMO(\mathbb R^n)}= \sup_{x \in \mathbb{R}^n,r>0} \fint_{B(x,r)} |f(s)-\{f\}_{B(x,r)}|\;ds
$$
and
$$
||f||_{BMO(\Omega)}=\sup_{x \in \mathbb{R}^n,r>0} 
\fint_{B(x,r)\cap \Omega} |f(s)-\{f\}_{B(x,r)\cap \Omega}|\;ds, $$
where $\{f\}_A = \fint_A f(s)\;ds=\frac{1}{|A|}\int_A ds$ are bounded.

By $C$ we denote a generic constant that is independent from $\nu$.

\section{Proofs of theorems}

We start with the proof of the estimate which plays the key role in next proofs.

\bigskip

\begin{proofoftheorem}{1.1} 
Consider $g \in \mathcal{H}^1(\mathbb{R}^n)$. By characterization of $\mathcal{H}^1(\mathbb{R}^n)$ we have
$$ \|g\|_{\mathcal{H}^1(\mathbb{R}^n)} = \|g\|_{L_1(\mathbb{R}^n)} + \sum_{k=1}^n\|R_k g\|_{L_1(\mathbb{R}^n)},$$
where the Riesz transform is given in the usual ways as $ \mathcal{F}[R_kf_k]=\frac{\xi_k}{|\xi|}\mathcal{F}[f_k]$. 
Using the fact that $BMO(\mathbb{R}^n)=(\mathcal{H}^1(\mathbb{R}^n))^*$ we get
\begin{equation}\label{estimate_2}
	\left| \int_{\mathbb{R}^n} fg \;dx \right| \leq \|f\|_{BMO(\mathbb{R}^n)}\|g\|_{\mathcal{H}^1(\mathbb{R}^n)}\leq \|f\|_{BMO(\mathbb{R}^n)}\left( \|g\|_{L_1(\mathbb{R}^n)} + \sum_{k=1}^n\|R_k g\|_{L_1(\mathbb{R}^n)} \right).
\end{equation}
For the characterization of the Hardy space $\mathcal{H}^1(\mathbb{R}^n)$ the reader may refer to \cite{ST1}. Hence it suffices to obtain an estimate on the $L_1$-norm of $R_k g$. We use the
classical Zygmund result that can be found in \cite{TOR1}.
\begin{proposition}
Let $h$ be a sufficiently smooth function with bounded
support. Then
\begin{equation}\label{estimate_Zygmund}
	\|R_k h\|_{L_1(\mathbb{R}^n)} \leq C+C\int_{\mathbb{R}^n} h\ln ^+ h\;dx
\end{equation}
where $\ln^+ a= \max\{ \ln a, 0 \}$ and the constant $C$ depends on the measure of support of $h$.
\end{proposition}

By elementary scaling we change inequality (\ref{estimate_Zygmund}) to get
$$ \|R_k g\|_{L_1(\mathbb{R}^n)} \leq \lambda +C\int_{\mathbb{R}^n} g \ln^+ g/\lambda \; dx$$
for any $\lambda \in \mathbb{R}_+$. Consider $\ln^+ g/\lambda = \ln g - \ln \lambda$ for $g \geq \lambda$, then
$$ \left|\ln (g|_{\{g \geq \lambda \}})\right| \leq \left|\ln(1+\|g\|_{L_\infty(\mathbb{R}^n)}) \right| + \left| \ln
\frac{g}{\|g\|_{L_\infty(\mathbb{R}^n)}+1} |_{\{g \geq \lambda\}}\right|.$$
Since $\frac{1}{1+\|g\|_{L_\infty(\mathbb{R}^n)}}\leq 1$ by elementary properties of logarithms we obtain
$$ \left| \ln(1+\|g\|_{L_\infty(\mathbb{R}^n)}) \right| + \left| \ln \frac{g}{\|g\|_{L_\infty(\mathbb{R}^n)}+1} |_{\{g \geq\lambda\}}\right| \leq \left| \ln(1+\|g\|_{L_\infty(\mathbb{R}^n)}) \right| + \left| \ln \frac{\lambda}{\|g\|_{L_\infty(\mathbb{R}^n)}+1} \right|.$$
Since it suffices to consider $\lambda \leq \|g\|_{L_\infty(\mathbb{R}^n)}$, we get
$$ |\ln (g|_{\{g \geq \lambda \}})| \leq 2 \ln(\|g\|_{L_\infty(\mathbb{R}^n)}+1)+|\ln \lambda|.$$
Choose $\lambda = \|g\|_{L_1(\mathbb{R}^n)}$. We then have
\begin{equation}\label{estimate_3}
	\|R_k g\|_{L_1(\mathbb{R}^n)} \leq c\|g\|_{L_1(\mathbb{R}^n)} + 2 \int_{\mathbb{R}^n} g \left(
\ln(\|g\|_{L_\infty(\mathbb{R}^n)}+1)+|\ln\|g\|_{L_1(\mathbb{R}^n)}|\right)\;dx
\end{equation}
$$ \leq c\|g\|_{L_1(\mathbb{R}^n)}\left( 1 + 2\ln(\|g\|_{L_\infty(\mathbb{R}^n)}+1) +
|\ln\|g\|_{L_1(\mathbb{R}^n)}|\right).$$
Inequality (\ref{ineq_est}) follows from inequalities (\ref{estimate_2}),(\ref{estimate_3}).
\end{proofoftheorem}

\bigskip

\begin{remark}
	Let $\Omega$ be a bounded subset of $\mathbb{R}^n$, $f \in BMO(\Omega)$, $g \in L_1(\Omega)\cap L_\infty(\Omega)$. Then 
		\begin{equation}\label{ineq_est_rem}
	\left|\int_{\Omega}fg\;dx \right|\leq C\|f\|_{BMO(\Omega)}\|g\|_{L_1(\Omega)}
	\left[  |\ln\|g\|_{L_1(\Omega)} |+\ln(1+\|g\|_{L_\infty(\Omega)})
	  \right].
\end{equation}
\end{remark}
\begin{proof}
	Extending $f,g$ by $0$ outside $\Omega$ we can apply Theorem 1.1. Now notice that for such 
	extension $\|g\|_{L_1(\mathbb{R}^n)}=\|g\|_{L_1(\Omega)}$, $\|g\|_{L_\infty(\mathbb{R}^n)}=\|g\|_{L_\infty(\Omega)}$ 
	and $\|f\|_{BMO(\mathbb{R}^n)} \leq C \|f\|_{BMO(\Omega)}$.
\end{proof}

\bigskip

Inequality (\ref{ineq_est}) is a key estimate in the proof of Theorem 1.2.

\bigskip

\begin{proofoftheorem}{1.2}
Let $\{u_1^E,p_1^E \}$ and $\{u_2^E,p_2^E \}$ be two different solutions of (\ref{Euler}). Subtracting the equations we get
\begin{equation}
                 (u_1^E-u_2^E)+ u_1^E\cdot
\nabla(u^E_1-u^E_2)+(u^E_1-u^E_2)\nabla u^E_2 = -\nabla(p_1^E-p_2^E).
         \end{equation}
         Multiplying both sides by $(u_1^E-u_2^E)$ and integrating over
$\Omega$ we obtain
         \begin{equation}\label{integration1}
                 \frac{1}{2}\frac{d}{dt}\int_\Omega (u^E_1-u^E_2)^2\;dx +
\int_\Omega(u^E_1-u^E_2)u^E_1\nabla(u^E_1-u^E_2)\;dx
         \end{equation}
         $$ + \int_\Omega(u^E_1-u^E_2)^2\nabla u^E_2\;dx = -
\int_\Omega(u^E_1-u^E_2)\nabla(p^E_1-p^E_2)\;dx.$$
         Integrating by parts, using boundary conditions and incompressibility
of flow we reduce (\ref{integration1}) to
         \begin{equation}\label{product1}
                 \frac{1}{2}\frac{d}{dt}\int_\Omega (u^E_1-u^E_2)^2\;dx +
\int_\Omega(u^E_1-u^E_2)^2\nabla u^E_2\;dx =0.
         \end{equation}
         Notice that from the assumptions on $u_1^E,u_2^E$ and (\ref{product1}) it follows that 
	 $\|u_1^E-u_2^E\|_{L_2(\Omega)} \in C([0,T])$. We split $\alpha= \alpha_m+\alpha_r$ where $|\alpha_m| = \min(|\alpha|,m)$
for some $m>1$. Notice that we can extend all functions outside $\Omega$ by $0$. Upon Theorem 1.1 and Remark 2.1 we get
         \begin{equation}\label{ineq1}
                 \left|\int_\Omega |\alpha\beta|\;dx \right| \leq C
\|\beta\|_{BMO(\Omega)}\|\alpha_m\|_{L_1(\Omega)}(1+|\ln
\|\alpha_m\|_{L_1(\Omega)}|+\ln(1+m))+\int_\Omega|\alpha_r\beta|\;dx.
         \end{equation}
	 Denote $f(t) = 2C\|\beta(t)\|_{BMO}$, $g(t)=\int_\Omega
|\alpha_r\beta|\;dx$, $x(t)=\|u^E_1-u^E_2\|^2_{L_2(\Omega)}$.
 Consider $x(t)$ small enough so that the function $|x(t)\ln x(t)|$ is increasing (which by the continuity of $x(t)$
 is equivalent to restricting our attention to sufficiently small $T$), then (\ref{ineq1}) can be restated as follows
 \begin{equation}\label{ineq2}
                 \left|\int_\Omega |\alpha\beta|\;dx \right| \leq C
\|\beta\|_{BMO(\Omega)}\|\alpha\|_{L_1(\Omega)}(1+|\ln
\|\alpha\|_{L_1(\Omega)}|+\ln(1+m))+\int_\Omega|\alpha_r\beta|\;dx.
         \end{equation}  
Thus from (\ref{ineq2}) we obtain the following inequality
\begin{equation}\label{diff_0}
								\begin{gathered}
                 \dot{x} \leq f(t)x(t)(|\ln x(t)|+1+\ln(1+m))+g(t), \hfill \\
                 x(0)=0. \hfill
                \end{gathered}
\end{equation}

To find a good estimate on $x(t)$ we  introduce the following equation
\begin{equation}\label{diff_1}
								\begin{gathered}
                 \dot{y} = f(t)y(t)(|\ln y(t)|+1+\ln(1+m))+g(t), \hfill \\
                 y(0)=1/m, \hfill
                \end{gathered}
\end{equation}
for some $m$ large enough. From the Osgood existence theorem we know that there exists a unique local solution to (\ref{diff_1}). 
Additianlly the r.h.s. of $(\ref{diff_1})_1$ guarantees that $y(\cdot)$ is increasing. It implies that
the solution of (\ref{diff_1}) majorizes $x(t)$, i.e.:
$$
0\leq x(t)\leq y(t) \mbox{ \ \ \ for \ \ } t\in [0,T].
$$ 

Hence we investiagate the bahavior of the solutions to  (\ref{diff_1}). By Gronwall's inequality we get
         \begin{equation}\label{ineq10}
         \begin{gathered}
y(t) \leq \frac{1}{m}\exp \left( \int_0^t f(s)\left( |\ln y(s)|+1+\ln(1+m)\right)\;ds \right) +  \\
\int_0^t g(s) \exp\left( \int_s^t f(\tau) \left( |\ln y(\tau)| + 1+\ln(1+m) \right) \;d\tau  \right)\;ds \leq \\
\frac{1}{m} \exp\left( (1+\ln(1+m)) \int_0^t f(s)\;ds \right)\exp\left( \int_0^t f(s)|\ln y(s)|\;ds\right)+ \\
+ \exp\left( \int_0^t f(\tau) |\ln y(\tau)|\;d\tau \right)\exp\left( (1+\ln(1+m)) \int_0^t f(\tau)\;d\tau \right)\int_0^t g(s)\;ds.
                 \end{gathered}
                 \end{equation}
Since $y(t)\geq 1/m$ implies $|\ln y(t)| \leq \ln m$ we can estimate the right hand-side (modulo some constant) of (\ref{ineq10}) by
                 \begin{equation}
                 \begin{gathered}
                 \frac{1}{m}(1+m)^{\int_0^t f(s)\;ds}m^{\int_0^t f(s)\;ds}+
m^{\int_0^t f(s)\;ds}(1+m)^{\int_0^t f(s)\;ds}\int_0^t g(s)\;ds = \\
                 (m(1+m))^{\int_0^t f(s)\;ds} \left( \frac{1}{m} +
\int_0^t g(s)\;ds \right) \leq \\
(2m^2)^{\int_0^t f(s)\;ds} \left( \frac{1}{m} +
\int_0^t \int_\Omega |\alpha_r \beta|\;dxds \right).
\end{gathered}
\end{equation}
This shows that it suffices to control the part $\alpha_r$. In this case we have estimates on the measure of the support.
Since $\alpha \in L_\infty(0,T;L_{1+\sigma/2}(\Omega))$ and $\beta \in BMO(\Omega)$ hence $\beta \in L_p(\Omega)$ for any $p<\infty$ we have by elementary H\"older's inequality
$$ \int_\Omega|\alpha_r \beta|\;dx \leq \|\alpha_r\|_{L_{1+\sigma/4}(\Omega)}\|\beta\|_{L_{(1+\sigma/4)'}(\Omega)},$$
hence we obtain a bound
\begin{equation}\label{estimate3}
\begin{gathered}
 (2m^2)^{\int_0^t f(s)\;ds} \left( \frac{1}{m} +
\int_0^t \int_\Omega |\alpha_r \beta|\;dxds \right)  \\
 \leq (2m^2)^{\int_0^t f(s)\;ds} \left( \frac{1}{m} +
\|\alpha_r\|_{L_\infty(0,T;L_{1+\sigma/4}(\Omega))}\|\beta\|_{L_1(0,T;L_{(1+\sigma/4)'}(\Omega))} \right).
\end{gathered}
\end{equation}
From Chebyschev inequality we have
\begin{equation}\label{cheb_1}
 |supp \; \alpha_r| \leq \left( \frac{\|\alpha\|_{L_\infty(0,T;L_{1+\sigma/2}(\Omega))}}{m}\right)^{1+\sigma/2}
\end{equation}
uniformly in time. Notice that by H\"olders inequality
\begin{equation}\label{hold_1}
	\|\alpha_r\|_{L_{1+\sigma/4}(\Omega)} = \left( \int_{supp\; \alpha_r} |\alpha_r|^{1+\sigma/4} \;dx\right)^{\frac{1}{1+\sigma/4}} \leq \left( |supp\;\alpha_r|^{\frac{\sigma}{4+2\sigma}} \;\cdot \; \|\alpha_r\|_{L_{1+\sigma/2}(\Omega)}^{\frac{4+\sigma}{4+2\sigma}} \right)^{\frac{1}{1+\sigma/4}}.
	\end{equation}
Inequalities (\ref{hold_1}) and (\ref{cheb_1}) imply
\begin{equation}
	\|\alpha_r\|_{L_{1+\sigma/4}(\Omega)} \leq m^{-\frac{\sigma}{4+\sigma}} \; \cdot \; \|\alpha_r\|_{L_{1+\sigma/2}(\Omega)}^{1-\frac{2\sigma}{(\sigma+2)(\sigma+4)}},
\end{equation}
hence
\begin{equation}
	\|\alpha_r\|_{L_\infty(0,T;L_{1+\sigma/4}(\Omega))} \leq m^{-\frac{\sigma}{4+\sigma}} \; \cdot \; \|\alpha_r\|_{L_\infty(0,T;L_{1+\sigma/2}(\Omega))}^{1-\frac{2\sigma}{(\sigma+2)(\sigma+4)}}.
\end{equation}
Choose $0<t_1\leq T$ small enough so that $4\int_0^{t_1}f(s)\;ds < \frac{\sigma}{4+\sigma}$, then for $0\leq t \leq t_1$ there is $y(t) \leq C m^{-\frac{\sigma}{8+2\sigma}}$.
Letting $m \to \infty$ we get 
$x(t)=0$ for $0\leq t \leq t_1$ which reads $u^E_1=u^E_2$ for $0\leq t \leq t_1$. We can continue this procedure 
starting at $t=t_1$ and get uniqueness for all $t\in[0,T]$. 
\end{proofoftheorem}

\bigskip

Estimate (\ref{ineq_est}) can also be used to give insight into the rate of convergence in the inviscid limit of the system (\ref{NS}).

\bigskip

\begin{proofoftheorem}{1.3}
	Let $\{u^\nu,p^\nu\}$ and $\{u^{E},p^{E}\}$ be solutions to problems (\ref{NS}) and (\ref{Euler}) respectively. Subtracting these equations we get
 \begin{equation}
                 (u^\nu-u^E)-\nu\Delta u^\nu+u^\nu\cdot
\nabla(u^\nu-u^E)+(u^\nu-u^E)\nabla u^E = -\nabla(p^\nu-p^E).
         \end{equation}
         Multiplying both sides by $(u^\nu-u^E)$ and integrating over
$\Omega$ we obtain
         \begin{equation}\label{integration2}
                 \frac{1}{2}\frac{d}{dt}\int_\Omega (u^\nu-u^E)^2\;dx + \nu \int_\Omega (u^\nu-u^E)\Delta u^\nu\;dx + 
\int_\Omega(u^\nu-u^E)u^\nu\nabla(u^\nu-u^E)\;dx
         \end{equation}
         $$ + \int_\Omega(u^\nu-u^E)^2\nabla u^E\;dx = -
\int_\Omega(u^\nu-u^E)\nabla(p^\nu-p^E)\;dx.$$
         Integrating by parts, using boundary conditions and incompressibility
of flow we reduce (\ref{integration2}) to
         \begin{equation}\label{product2}
                 \frac{1}{2}\frac{d}{dt}\int_\Omega (u^\nu-u^E)^2\;dx +
\int_\Omega(u^\nu-u^E)^2\nabla u^E\;dx - \nu \int_\Omega \nabla u^\nu \nabla (u^\nu-u^E)\;dx =0.
         \end{equation}
Let $\alpha=(u^\nu-u^E)^2$, $\beta= \nabla u^E$. 
Extending all functions by $0$ outside $\Omega$, by Theorem 1.1 and Remark 2.1 we have
 \begin{equation}\label{rate_ineq_2}
                 \left|\int_\Omega |\alpha\beta|\;dx \right| \leq C
\|\beta\|_{BMO}\|\alpha\|_{L_1(\Omega)}(1+|\ln
\|\alpha\|_{L_1(\Omega)}|+\ln(1+\|\alpha\|_{L_\infty}(\Omega))).
         \end{equation}
Let $\alpha=\alpha_\nu+\alpha_r$ so that $|\alpha_\nu| = \min(\frac{1}{\nu},|\alpha|)$. Proceeding as in the proof of Theorem 1.2 we denote $f(t) = 2C\|\beta(t)\|_{BMO(\Omega)}$, $g(t)=\int_\Omega
|\alpha_r\beta|\;dx$, $x(t)=\|u^\nu-u^E\|^2_{L_2(\Omega)}$. 
From (\ref{rate_ineq_2}) we get the following inequality
\begin{equation}\label{diff_20}
\begin{gathered}
	             \dot{x} \leq f(t)x(t)\left(|\ln x(t)|+1+\ln\left(1+\frac{1}{\nu}\right)\right)+g(t) + \nu g_0(t)^2, \hfill \\
                 x(0)=0. \hfill
\end{gathered}
\end{equation}

To find a good estimate on $x(t)$ we  introduce the following equation
\begin{equation}\label{diff_21}
								\begin{gathered}
                 \dot{y} = f(t)y(t)\left(|\ln y(t)|+1+\ln\left(1+\frac{1}{\nu}\right)\right)+g(t) + \nu g_0(t)^2, \hfill \\
                 y(0)=\nu \hfill
                \end{gathered}
\end{equation}
for some $\nu$ sufficiently small. From the Osgood existence theorem we know that there exists a unique local solution to (\ref{diff_21}). 
The solution of (\ref{diff_21}) majorizes $x(t)$, i.e.:
$$
x(t)\leq y(t) \mbox{ \ \ \ for \ \ } t\in [0,T].
$$ 
From (\ref{diff_21}) we have by Gronwall's inequality
\begin{equation}\label{2.24}
\begin{gathered}
y(t) \leq \nu\exp \left( \int_0^t f(s)\left[ |\ln y(s)|+1+\ln\left(1+\frac{1}{\nu}\right)\right]\;ds \right) + \\
\int_0^t \left(g(s)+g_0(s)^2\nu\right) \exp\left( \int_s^t f(\tau) \left[ |\ln y(\tau)| + 1+\ln\left(1+\frac{1}{\nu}\right) \;d\tau \right] \right)\;ds \leq \\
 \nu \exp\left( \left(1+\ln\left(1+\frac{1}{\nu}\right)\right) \int_0^t f(s)\;ds \right)\exp\left( \int_0^t f(s)|\ln y(s)|\;ds\right)+  \\
 \exp\left( \left(1+\ln\left(1+\frac{1}{\nu}\right)\right)\int_0^t f(\tau)\;d\tau \right)\exp\left( \int_0^t f(s) |\ln y(s)|\;ds \right) \\ \times \int_0^t (g(s)+g_0(s)^2\nu)\;ds.
\end{gathered}
\end{equation}
The condition $y(t) \geq \nu$ for sufficiently small $\nu$ gives $|\ln y(t)| \leq -\ln \nu = \ln \frac{1}{\nu}$. Also let $\nu$ be small enough so that $\frac{1}{\nu}\left(1+\frac{1}{\nu}\right) \leq \frac{2}{\nu^2}$, thus we can estimate the right hand-side of (\ref{2.24}) (modulo a constant) by
\begin{equation}
\begin{gathered}
 \nu\left(\frac{2}{\nu^2}\right)^{\int_0^t f(s)\;ds}+
\left(\frac{2}{\nu^2}\right)^{\int_0^t f(s)\;ds}\int_0^t (g(s)+\nu g_0(s)^2)\;ds = \\   
= \left(\frac{2}{\nu^2}\right)^{\int_0^t f(s)\;ds} \left( \nu +
\int_0^t \int_\Omega |\alpha_r \beta|\;dxds  + \nu \int_0^t g_0(s)^2\;ds \;dt \right).
\end{gathered}
\end{equation}
Since $\alpha \in L_\infty(0,T;L_{1+\sigma/2}(\Omega))$ and $\beta \in BMO(\Omega)$ hence $\beta \in L_p(\Omega)$ for any $p<\infty$ we have 
$$ \int_\Omega|\alpha_r \beta|\;dx \leq \|\alpha_r\|_{L_{1+\sigma/4}(\Omega)}\|\beta\|_{L_{(1+\sigma/4)'}(\Omega)},$$
thus
\begin{equation}\label{estimate31}
\begin{gathered}
\left(\frac{2}{\nu^2}\right)^{\int_0^t f(s)\;ds} \left( \nu +
\int_0^t \int_\Omega |\alpha_r \beta|\;dxds + \nu \int_0^t |g_0(s)|^2  \;ds \right) \leq \\
 \left(\frac{2}{\nu^2}\right)^{\int_0^t f(s)\;ds} \left( \nu +
\|\alpha_r\|_{L_\infty(0,T;L_{1+\sigma/4}(\Omega))}\|\beta\|_{L_1(0,T;L_{(1+\sigma/4)'}(\Omega))} + \nu\|g_0\|^2_{L_2(0,T)}\right).
\end{gathered}
\end{equation}
From the Chebyshev inequality we notice that
\begin{equation}\label{cheb_2}
 |supp \; \alpha_r| \leq \left( \nu \|\alpha\|_{L_\infty(0,T;L_{1+\sigma/2}(\Omega))} \right)^{1+\sigma/2},
\end{equation}
uniformly in time. Notice that by H\"olders inequality
\begin{equation}\label{hold_2}
	\|\alpha_r\|_{L_{1+\sigma/4}(\Omega)} = 
	\left( \int_{supp\; \alpha_r} |\alpha_r|^{1+\sigma/4} \;dx\right)^{\frac{1}{1+\sigma/4}}
	 \leq \left( |supp\;\alpha_r|^{\frac{\sigma}{4+2\sigma}} \;\cdot \;
	  \|\alpha_r\|_{L_{1+\sigma/2}(\Omega)}^{\frac{4+\sigma}{4+2\sigma}} \right)^{\frac{1}{1+\sigma/4}}.
	\end{equation}
Inequalities (\ref{hold_2}) and (\ref{cheb_2}) imply
\begin{equation}
	\|\alpha_r\|_{L_{1+\sigma/4}(\Omega)} \leq \nu^{\frac{\sigma}{4+\sigma}} \; 
	\cdot \; \|\alpha_r\|_{L_{1+\sigma/2}(\Omega)}^{1-\frac{2\sigma}{(\sigma+2)(\sigma+4)}},
\end{equation}
hence
\begin{equation}
	\|\alpha_r\|_{L_\infty(0,T;L_{1+\sigma/4}(\Omega))} \leq \nu^{\frac{\sigma}{4+\sigma}} \; 
	\cdot \; \|\alpha_r\|_{L_\infty(0,T;L_{1+\sigma/2}(\Omega))}^{1-\frac{2\sigma}{(\sigma+2)(\sigma+4)}}.
\end{equation}
Choose $0<t_1\leq T$ small enough so that $4\int_0^{t_1}f(s)\;ds < \frac{\sigma}{4+\sigma}$, then for 
$0\leq t \leq t_1$ there is $y(t) \leq C \nu^{\frac{\sigma}{8+2\sigma}}$. Consider now, $t_1 \leq t \leq T$ and 
a problem analogous to (\ref{diff_21}) but with initial condition $y(t_1)=\nu^{\frac{\sigma}{8+2\sigma}}$. 
Repeating all above estimates we pick $t_1 < t_2 \leq T$ such that 
$\sup_{t_1 \leq t \leq t_2}\|u^\nu-u^E\|_{L_2(\Omega)}\leq C \nu^{\frac{\sigma}{16+4\sigma}}$. Due to integrability of $f(t)$ iterating the 
procedure we eventually cover the whole interval $[0,T]$. This way we obtain the explicit rate of the convergence which
depends mainly on the structure of integrability of function $f_0$. Thus we proved (\ref{1.9}).

Assuming additionally condition (\ref{1.10}) improves the result and gives an explicit uniform rate of convergence. 
The basic estimate presented above gives 
\begin{equation}
	x(t) \leq C\nu\left(\frac{2}{\nu^2}\right)^{Mt}
\end{equation}
with $M$ as in (\ref{1.10}).
Fix some $n \in \mathbb{N}$ and consider $0 \leq t \leq T/(2Mn)$. We have
 $\sup_{0\leq t \leq T/(2Mn)} \|u^\nu-u^E\|_{L_2(\Omega)} \leq C \nu^{1-T/n}$. The time interval has been divided into 
 $2Mn$ parts and repeating the estimate we get
\begin{equation}
	\sup_{0\leq t \leq T} \|u^\nu-u^E\|_{L_2(\Omega)} \leq C \nu^{(1-T/n)^{2Mn}},
\end{equation} 
and taking limit $n \to \infty$ we obtain 
\begin{equation}
	\sup_{0\leq t \leq T} \|u^\nu-u^E\|_{L_2(\Omega)} \leq C \nu^{e^{-2MT}}.
\end{equation} 
Theorem 1.3 is proved.

\end{proofoftheorem}

\noindent{\footnotesize {\bf Acknowledgments.} The first author (PBM) has been partly supported by
Polish KBN grant No. 1 P03A 021 30.}


{\footnotesize
\begin{thebibliography}{99}

\bibitem {Cl} Clopeau T., Mikeli\'c A., Robert R., \textit{On the vanishing viscosity limit for the 2D incompressible 
			       Navier--Stokes equations with friction type boundary conditions}. Nonlinearity {\bf 11} (1998), 1625--1636.

\bibitem {CONST} Constantin P., Wu J., \textit{Inviscid limit for vortex patches}. Nonlinearity {\bf 8} (1995), 735--742.

\bibitem {ST1} Fefferman C., Stein E.M., \textit{$H\sp{p}$ spaces of several variables}. Acta Math. {\bf 129} no. 3-4 (1972), 137--193.

\bibitem {HART} Hartman F., \textit{Ordinary differential equations}. John Wiley \& Sons, NY-London-Sydney (1964).

\bibitem {KAT1} Kato T., \textit{On classical solutions of the two-dimensional nonstationary Euler equation}. Arch. Rational Mech. Anal. {\bf 25} (1967), 188--200. 	


\bibitem {KOZ} Kozono H., Taniuchi Y., \textit{Limiting case of the Sobolev inequality in BMO, with application to the Euler equations}.  Comm. Math. Phys.  {\bf 214} no. 1 (2000), 191--200.

\bibitem {MU} Mucha P.B.,\textit{On the inviscid limit of the Navier-Stokes
equations for flows with large flux}.  Nonlinearity  {\bf 16} no. 5  (2003)
1715--1732.
								
\bibitem {Ru} Rusin W.M.,\textit{On the inviscid limit for the solutions of two-dimensional incompressible Navier-Stokes equations with slip-type boundary conditions}. Nonlinearity {\bf 19} no. 6 (2006), 1349--1363. 
				
\bibitem {TOR1} Torchinsky A., \textit{Real-variable methods in harmonic analysis}. Pure and Applied Mathematics, 123. Academic Press, Inc., Orlando, FL, 1986.				
								
\bibitem {Y1} Yudovich V., \textit{Nonstationary flow of an ideal incompressible liquid}. Zhurn. Vych. Mat. {\bf 3} (1963), 1032--1066.
				
\bibitem {Y2} Yudovich V., \textit{Uniqueness theorem for the basic nonstationary problem in the dynamics of an ideal incompressible fluid}. Math. Res. Lett. {\bf 2} (1995), 27--38. 
				
\bibitem {Vish} Vishik M., \textit{Incompressible flows of an ideal fluid with vorticity in borderline spaces of Besov type}.	Ann. Sci. cole Norm. Sup. (4) {\bf 32} no. 6 (1999), 769--812. 
				
\bibitem {Zyg1} Zygmund A., \textit{Trygonometric Series}. Cambridge Univ. Press, London-NY, 1959.
				
				
\end{thebibliography}
}
\end{document}